\documentclass{amsart}[12pt,a4paper]

\usepackage{amssymb,amsmath,amsthm, amscd, enumerate, mathrsfs}
\usepackage{graphicx, hhline}
\usepackage[colorlinks=true,pagebackref,hyperindex]{hyperref}
\usepackage[all]{xy}
\usepackage{color}
\usepackage[backrefs, alphabetic, initials]{amsrefs}
\usepackage{color}  
\usepackage{latexsym}
\usepackage[T1]{fontenc} 
\usepackage{fancyhdr}
\usepackage{enumerate}

\title[Rational curves and strictly nef divisors on Calabi--Yau 3-folds]
{Rational curves and strictly nef divisors on Calabi--Yau threefolds}
\author{Haidong Liu and Roberto Svaldi}

\date{2020/10/17, version 0.03}



\address{Peking University, Beijing International Center for Mathematical Research, 
Beijing, 100871, China}
\email{hdliu@bicmr.pku.edu.cn}

\address{EPFL, SB MATH-GE, MA B1 497 (Bâtiment MA), Station 8, CH-1015 Lausanne, Switzerland}
\email{roberto.svaldi@epfl.ch}

\DeclareMathOperator{\nefeff}{Nef^{e}}

\DeclareMathOperator{\NS}{N^1}

\DeclareMathOperator{\Eff}{\overline{Eff}}

\DeclareMathOperator{\NE}{\overline{NE}}
\DeclareMathOperator{\eff}{Eff}

\DeclareMathOperator{\Aut}{Aut}

\DeclareMathOperator{\nef}{Nef}

\DeclareMathOperator{\fol}{\mathcal F}
\DeclareMathOperator{\strutt}{\mathcal{O}}

\newtheorem{thm}{Theorem}[section]
\newtheorem{lem}[thm]{Lemma}
\newtheorem{prop}[thm]{Proposition}
\newtheorem{conj}[thm]{Conjecture}
\newtheorem{cor}[thm]{Corollary}

\theoremstyle{definition}

\newtheorem{rem}[thm]{Remark}
\newtheorem*{ack}{Acknowledgments}

\newtheorem*{claim}{Claim}

\makeatletter
    
    \@addtoreset{equation}{section}
\makeatother

\begin{document}

\begin{abstract}
We give a criterion for a nef divisor $D$ to be semiample on a Calabi--Yau threefold $X$ when $D^3=0=c_2(X)\cdot D$ and $c_3(X)\neq 0$. 
As a direct consequence, we show that on such a variety $X$, if $D$ is strictly nef and $\nu(D)\neq 1$, then $D$ is ample;
we also show that if there exists a nef non-ample divisor $D$ with $D\not\equiv 0$, then $X$ contains a rational curve when its topological Euler characteristic is not $0$.
\end{abstract}

\thanks{
Part of this work was carried out while RS was visiting Professor Jacopo Stoppa at SISSA, Trieste. 
He would like to thank Professor Stoppa and SISSA for providing such a congenial place to work. 
The visit was supported by funding from the European Research Council under the European Union's Seventh Framework Programme (FP7/2007-2013)/ERC Grant agreement no. 307119.
The author was partially supported by Churchill College, Cambridge and from the European Union's Horizon 2020 research and innovation programme under the Marie Sk\l{}odowska-Curie grant agreement No. 842071.
}

\maketitle

\setcounter{tocdepth}{1}
\tableofcontents

\section{Introduction}
\label{sec1}

From the point of view of the birational classification of algebraic varieties, varieties with trivial canonical bundle and mild singularities (in short, K-trivial varieties) are one of the fundamental building blocks. Hence, their study is central to the understanding of the structure of algebraic varieties as a whole. 
By a well-known result of Beauville and Bogomolov~\cite{beauville}, every smooth variety with trivial canonical bundle can be decomposed, after a finite \'etale cover, as a product of abelian, hyperk\"ahler or Calabi--Yau manifolds. 
For the purpose of this note, a smooth projective variety $X$ with trivial canonical bundle is Calabi--Yau if it is simply connected and $H^i(X, \mathcal{O}_X)=0$ for $0<i<\dim X$.

One of the central problems in modern birational geometry is the so-called Abundance Conjecture, see~\cite{kollar.sing}*{\S 8.7}.

\begin{conj}
\label{conj.ab}
Let $(X, \Delta)$ be a log canonical pair, where $\Delta$ is a $\mathbb Q$-divisor. 
If $K_X+\Delta$ is nef, then it is semiample.
\end{conj}

For K-trivial varieties, the above conjecture is expected to hold in even greater generality, cf.~\cite{kollar}*{Conjecture 51},~\cite{lop1}*{Section 4} and ~\cite{lp18} for conventions and a survey of several related conjectures.

\begin{conj}
[Semiampleness conjecture for K-trivial varieties]
\label{conj-abundance}
Let $X$ be a projective klt variety with $K_X \equiv 0$. Assume that $H^1(X, \mathcal{O}_X)=0$. Let $D$ be a nef divisor on $X$. Then $D$ is semiample.
\end{conj}

The major difference between the Abundance Conjecture~\ref{conj.ab} for a K-trivial variety and the Semiampleness Conjecture~\ref{conj-abundance} is that the divisor $D$, while nef, is not required to be effective. 
The assumption that $H^1(X, \mathcal{O}_X)=0$ implies that ${\rm Pic}^0(X)$ is $0$-dimensional, which is a necessary condition, as otherwise there would be numerically trivial line bundles that are non-torsion. The condition is automatically satisfied for Calabi--Yau manifolds and simple hyperk\"ahler manifolds. 
In the realm of simple hyperk\"ahler manifolds, Conjecture~\ref{conj-abundance} is a kind of the so-called Strominger--Yau--Zaslow (SYZ) conjecture, see~\cite{ver}*{Conjecture~1.7}.

While the Abundance Conjecture is known to hold in a number of cases, e.g., it holds for projective K-trivial threefolds, so that if $D$ is an effective nef divisor then $D$ is semiample, cf.~\cite{kmm}*{Corollary, p.100} or~\cite{oguiso}, Conjecture~\ref{conj-abundance} is only known to hold in full generality for K-trivial varieties of dimension at most 2.
In dimension 3 and higher, only very few cases of the conjecture have been verified, cf.~\cites{kollar, lop1}.

In this paper, we study some cases of Conjecture~\ref{conj-abundance} together with some applications to other problems on Calabi--Yau threefolds.

The first result that we present is the following theorem quantifying precisely the possible failure of the Semiampleness Conjecture for a nef divisor of numerical dimension 2 on a Calabi--Yau threefold.
We recall the definitions for the readers that 
a prime divisor $S$ is said to be {\it orthogonal} to a Cartier divisor $D$ if $S \cdot D=0$;
the {\em{numerical dimension}} $\nu(D)$ of a nef divisor $D$ is $\nu(D):=\max\{h\in \mathbb N \ \vert \ D^h\not\equiv 0\}$.

\begin{thm}[Theorem~\ref{lop.thm1}]
\label{lop.thm}
Let $X$ be a Calabi--Yau threefold and $D$ a nef divisor on $X$ with $\nu(D)=2$. Let $S_1, \ldots, S_r$ be all the prime divisors on $X$ orthogonal to $D$. Let $g_j$ be the irregularity of a resolution of $S_j$. If
\[
\frac{c_3(X)}{2} \neq r-\sum_{j=1}^r g_j,
\]
then $D$ is semiample.
\end{thm}

Theorem~\ref{lop.thm} is a generalization of \cite{lop}*{Theorem~8.5}.
The important improvement in this generalization is that we have dropped the conditions~\cite{lop}*{Theorem~8.5(i)-(ii)} and obtained an unconditional statement.
The ideas and statements that appear here and in~\cite{lop} are influenced by early work of Wilson on this topic, cf.~\cites{wilson2, wilson}.

Conjecture~\ref{conj-abundance} is strictly intertwined with the existence of rational curves on K-trivial varieties, a classically hard problem particularly on Calabi--Yau varieties. This link is highlighted by the following conjecture of Oguiso.

\begin{conj}
[\cite{oguiso}*{Conjecture, p.456}]
\label{og.rat.curv.conj.intro}
	Let $X$ be a Calabi--Yau threefold. 
	Assume that there exists a non-trivial Cartier divisor 
	$D$ contained in the boundary of $\nef(X)$.
	Then $X$ contains a rational curve.
\end{conj}

Indeed, it is not hard to show that if Conjecture~\ref{conj-abundance} holds then also Oguiso's conjecture holds, as the Iitaka fibration of the divisor $D$ will contain rational curves in some fibers, cf.~Remark~\ref{rem.rc.fib}.
By applying the techniques used for Theorem~\ref{lop.thm}, we are able to prove the following almost complete solution of Conjecture~\ref{og.rat.curv.conj.intro}, which also complements the main result of~\cite{df}.

\begin{thm}
[{Theorem~\ref{rc-thm}}]
\label{rc-cor}
Let $X$ be a Calabi--Yau threefold. 
Assume that there exists a non-trivial Cartier divisor $D$ contained in the boundary of $\nef(X)$. 
Then $X$ contains a rational curve provided that either $c_3(X) \neq 0$ or the second Betti number $h^2(X) \neq 2,3,4$.
\end{thm}

In Section~\ref{sec7}, we show also how the problem of the existence of rational curves on Calabi--Yau threefolds is strictly intertwined with another important conjecture regarding the birational structure and the mirror symmetry of Calabi--Yau threefolds, the Kawamata--Morrison Conjecture. 
In this regard, we discuss also how geometric and cohomological considerations inspired by the conjecture could be used towards a definitive solution of Oguiso's conjecture.

Conjecture~\ref{conj-abundance} is also intertwined with another conjecture due to Serrano~\cite{serrano}. A {\it strictly nef} divisor $D$ on a normal variety $X$ is a nef Cartier divisor such that $D\cdot C>0$ for any curve $C$ on $X$. In~\cite{serrano}, Serrano conjectured that by twisting a sufficiently large multiple of a strictly nef divisor with the canonical bundle, one should obtain an ample divisor. 
On a Calabi--Yau manifold, as the canonical bundle is trivial, Serrano's conjecture just predicts the lack of non-ample strictly nef divisors, in agreement with the Semiampleness Conjecture~\ref{conj-abundance}.
Here we show an instance in which this principle holds on Calabi--Yau threefolds.

\begin{thm}[Theorem~\ref{sn-thm}]\label{sn-cor}
Let $X$ be a Calabi--Yau threefold with $c_3(X)\neq 0$. 
Let $D$ be a strictly nef divisor on $X$ with $\nu(D)\neq 1$.
Then $D$ is ample.
\end{thm}

The two results above are based on the solution of the following special case of Conjecture~\ref{conj-abundance} for Calabi--Yau threefolds, which can be seen as a corollary of Theorem~\ref{lop.thm}.

\begin{thm}\label{main.thm}
Let $X$ be a Calabi--Yau threefold with $c_3(X)\neq 0$. 
Let $D$ be a nef divisor with $D^3 = 0 = c_2(X) \cdot D$. 
Assume that there exists a very ample line bundle $H$ on $X$ and a general member $G \in |H|$ of the linear system such that $D|_G$ is ample.
Then $D$ is semiample.
\end{thm}

By the hypothesis of Theorem~\ref{main.thm}, the numerical dimension of $D$ is 2, that is, $\nu(D)=2$.
A weaker statement had been obtained in~\cite{lop}*{Theorem~1.2}: in their statement, the authors used a strong positivity condition on a certain sheaf of 1-forms with poles along an ample divisor and coefficients in $D$.
Hence, we successfully remove their hypothesis from the statement to obtain a more general result that allows us to obtain Theorem~\ref{rc-cor} and~\ref{sn-cor}.

The main idea in the proof of Theorem~\ref{main.thm} is to use the positivity condition on the restriction of $D$ to an element $G \in \vert H \vert$  together with Fujita's Vanishing Theorem~\ref{fujita-van} and the Hirzebruch--Riemann--Roch formula to produce sections of the nef divisor $D$.

Finally, if a Calabi--Yau threefold $X$ satisfies $c_3(X)=e(X)=2(h^{1,1}(X)-h^{2,1}(X))=0$, then such a $X$ would lie exactly on the axis of symmetry of the famous ``Zoo of Calabi--Yau threefolds'', cf.~\cite{ks}*{Figure 1}.
It is a folklore conjecture, based on the knowledge of existing models acquired over the past 40 years, that all Hodge numbers of Calabi--Yau threefolds should lie inside the cup shape appearing in~\cite{ks}*{Figure 1}.
One of the first sources for this cornucopia of examples to be explored were smooth divisors in the anti-canonical linear system of a 4-dimensional $\mathbb Q$-factorial Gorenstein toric Fano variety with terminal singularities.
By Batyrev's construction in~\cite{batyrev}, this class of toric Fano varieties are completely determined by reflexive polytopes in dimension 4; Kreuzer and Skarke in~\cite{ks} listed all elements in this class of polytopes, by means of a computer-aided classification.
As a result, we now know that there are $30,108$ distinct Hodge diamonds for the $473.8$ million families of Calabi--Yau threefolds.
Appearing in this list, there are just $136$ pairs of Hodge numbers of Calabi--Yau threefolds $X$ with $h^{1,1}(X)=h^{2,1}(X)$, hence, $c_3(X)=0$.
Thus, the condition that $c_3(X)=0$ appears relatively rarely distributed within our current picture of the ``Zoo of Calabi--Yau threefolds''.

\begin{ack}
The authors wish to thank Vladimir Lazi\'c for encouragement and useful conversations, and for providing a lot of suggestions after reading a first draft of this paper.
HL would like to thank Chen Jiang for helpful discussion and suggestions. 
Part of this project started as an appendix to~\cite{dfm}. 
RS is very grateful to Simone Diverio, Claudio Fontanari, and Diletta Martinelli for their encouragement and support as well as useful discussion and suggestions on the structure and the content of the material here included.
\end{ack}

\section{Preliminaries}\label{sec2}
Throughout this paper,  we work over the complex number field $\mathbb C$.
A {\em scheme} is always assumed to be separated and of finite type over $\mathbb{C}$. A {\em variety} is a reduced and irreducible scheme. 
For notation and conventions, we refer to~\cite{kollar-mori} and ~\cites{laz, laz2}.

As mentioned in the Introduction, we define the {\em{numerical dimension}} $\nu(D)$ of a nef $\mathbb Q$-divisor $D$ as
\[
\nu(D):=\max\{h\in \mathbb N \; |  \;D^h\not\equiv 0\}.
\]
It is well-known that $\nu(D) \geq 0$, $\kappa (D)\leq \nu(D) \leq \dim X$ and $\nu(D)=\dim X$ if and only if $D$ is big.
By Kleiman's Ampleness Criterion, it is also well-known that
\begin{lem}\label{positive.lemma}
If $\nu(D)=\dim X-1$, then 
$D^{\dim X-1}\cdot H>0$ for any ample divisor $H$ on $X$.
\end{lem}

\subsection{Cones of divisors}
We adopt the following notation to denote the various cones of divisors on a variety $X$:
\begin{itemize}
    \item $\NS(X)_{\mathbb R}$ denotes the N\'eron--Severi $\mathbb R$-vector space on $X$, while $\mathrm{N}_1(X)_{\mathbb R}$ denotes its dual, that is, the $\mathbb R$-vector space of 1-cycle classes modulo numerical equivalences.

    \item The {\em{effective cone}} $\eff(X)\subset \NS(X)_{\mathbb R}$ is the convex cone spanned by the classes of all effective $\mathbb R$-divisor classes on $X$.

    \item The {\em{pseudoeffective cone}} $\Eff(X)\subset \NS(X)_{\mathbb R}$ is the closure of $\eff(X)$. 
    The interior of $\Eff(X)$ is the convex cone of all big $\mathbb R$-divisor classes. 

    \item The {\em{nef cone}} $\nef(X)\subset \NS(X)_{\mathbb R}$ is the convex cone of all nef $\mathbb R$-divisor classes on $X$. 
    The interior of $\nef(X)$ is the convex cone of all ample $\mathbb R$-divisor classes. \\
    The {\em{effective nef cone}} $\nefeff(X)$ is equal to $\nef(X) \cap \eff(X)$. 
    It is not necessarily closed.
    
    \item The {\em cone of effective curves} $\NE(X) \subset \mathrm{N}_1(X)_{\mathbb R}$ is the closure of the cone spanned by effective 1-cycles with coefficients in $\mathbb R$.
\end{itemize}
It is well-known that $\nef(X)\subset \Eff(X)$ and $\nef(X)$ is the dual cone of $\NE(X)$, see~\cite{laz}*{Chapter~1}.

\subsection{Calabi--Yau threefolds}
\label{cy.sect}

We collect here a few basic results on Calabi--Yau threefolds that will be used in the rest of the paper, cf,~\cite{lop}. 
In this subsection, $X$ will denote a Calabi--Yau threefold.

In this paper, we will use the Hirzebruch--Riemann--Roch formula on a Calabi--Yau threefold $X$.
If $\mathcal F$ is a vector bundle on $X$, then the Hirzebruch--Riemann--Roch theorem formula is as follows
\begin{align}
\label{RR.cy.eq.vect}
\chi(X,\mathcal F)=\frac{1}{12}c_1(\mathcal{F})c_2(X)+\frac16\big(c_1(\mathcal{F})^3-3c_1(\mathcal{F})c_2(\mathcal{F})+3c_3(\mathcal{F})\big).
\end{align}
In particular, for a Cartier divisor $D$ on $X$,~\eqref{RR.cy.eq.vect} can be simplified to give
\begin{align}
\label{RR.cy.eq.div}
\chi(X,D)=\frac{D^3}{6}+\frac{D\cdot c_2(X)}{12}.
\end{align}

It is well-known that on $X$ the second Chern class $c_2(X)$ has non-negative intersection with nef divisors by the Miyaoka--Yau inequality~\cite{miyaoka}; moreover, $c_2(X) \neq 0$, as otherwise $X$ would be dominated by an abelian threefold, cf.~\cite{MR3522654}. 

Let us now consider a nef Cartier divisor $D$ on $X$.
By the solution of the Abundance Conjecture~\ref{conj.ab} for threefolds~\cites{kmm, kmm.err}, we know that $D$ is semiample on $X$ if and only if $\kappa(D) \geq 0$.
Therefore, the only case left to consider for Conjecture~\ref{conj-abundance} is when $\kappa(D)=-\infty$.
As $D$ is nef, this immediately implies that $D^3 =0$.
Using Kawamata--Viehweg vanishing and its generalization~\cite{laz}*{Example~4.3.7},  
if $\nu(D)=2$ and $c_2(X)\cdot D>0$, then $\kappa(D)\geq 0$ by~\eqref{RR.cy.eq.div} and $D$ is semiample.
Analogously, Oguiso showed that $c_2(X) \cdot D > 0$ implies that $D$ is semiample also when $\nu(D)=1$, see~\cite{oguiso}*{Main Theorem}.

In the proof of Theorems~\ref{main.thm} and~\ref{lop.thm}, we will use the following result taken from~\cite{wilson}*{Proposition 3.1},~\cite{lop}*{Proposition 3.4}.
A generalization of this result to higher dimensions appears in~\cite{lp}*{Theorem~8.1}.

\begin{prop}
\label{LP.van.prop}
 Let $X$ be a Calabi--Yau threefold. 
 Let $D$ be a nef divisor on $X$ with $\kappa(D) = -\infty$. 
 Then there exists a positive integer $m_0 \in \mathbb{N}_{>0}$ such that 
 \[
 H^0(X, \Omega_X^q \otimes \strutt_X(mD)) = 0
 \]
for all $\vert m \vert\geq m_0$ and $q\geq 0$.
\end{prop}

The following observation will play an important role in our discussion of rational curves on a Calabi--Yau varieties, cf. in particular~\S~\ref{sec6}.
\begin{rem}
\label{rem.rc.fib}
Let $D$ be a nef and semiample divisor on a Calabi--Yau threefold $X$, e.g., $D^3 > 0$ or $c_2(X) \cdot D >0$.
Then $X$ contains a rational curve.
To see this, it suffices to consider the Iitaka fibration $\phi \colon X \to Y$ induced by $|mD|$ for $m \gg 0$, since $D$ is semiample.
Then $\phi$ is either 
\begin{itemize}
    \item 
    a birational morphism,
    
    \item
    an elliptic fibration, or
    
    \item 
    a fibration in abelian/$K3$ surfaces.
\end{itemize} 
If $\phi$ is birational, then the exceptional locus of the birational modification is known to be uniruled.
For a fibration in $K3$ surfaces, it is well-known that the general fiber will contain rational curves.
In the other two cases, the results in the appendix of~\cite{oguiso} imply the existence of rational curves on $X$ at once, cf. also~\cite{dfm}.
\end{rem}

\subsection{Nef divisors and their null locus}

Let $X$ be a normal projective variety and  $D$ a nef Cartier divisor.
We study the set
$\mathcal{C}_D:=\{C \subset X \; | \; D \cdot C=0, \; \dim C =1, \; C \; \textrm{is irreducible}\}$.
We recall here the following standard general observation.
\begin{lem}\label{curve.lemma}
Let $X$ be a normal projective variety of dimension $3$ and $D$ a nef Cartier divisor on $X$ with $\nu(D)=2$.
Assume that the set $\mathcal{C}_D$ is at most countable.
Then for any very ample divisor $H$ on $X$ and for a very general member $G \in |H|$, $D|_G$ is ample.
\end{lem}

\begin{proof}
The proof of~\cite{lop}*{Lemma 2.5(b)} applies in this context.

\end{proof}

\begin{cor}\label{exist-cor}
Let $X$ be a normal projective variety of dimension $3$ and $D$ a strictly nef divisor on $X$ with $\nu(D)=2$.
Then for any very ample divisor $H$ on $X$ and for any very general member $G \in |H|$, $D|_G$ is ample.
\end{cor}

\begin{proof}
The result follows immediately from Lemma~\ref{curve.lemma} as $\mathcal{C}_D=\emptyset$ since $D$ is strictly nef.
\end{proof}

We will also use, particularly in \S~\ref{sec6}, the following result inspired by~\cite{lop}*{Theorem~2.7}.

\begin{lem}\label{h.sch.lemma}
Let $X$ be a Calabi--Yau threefold and $D$ a nef divisor on $X$ with $\nu(D)=2$. Assume that $X$ does not contain any rational curve. Then exactly one of the following conditions holds true:
  \begin{enumerate}
\item[1)] through any point $x \in X$ there exists a curve $C_x$ with $D \cdot C_x = 0$; in particular in this case $D$ is semiample; or
\item[2)] the set $\mathcal{C}_D$ is at most countable.
 \end{enumerate}
\end{lem}	

\begin{proof}
By the existence and properness of the Hilbert scheme, cf.~\cite{kollar.rational}*{Chapter~I}, and the nefness of $D$, if through a general point $x\in X$ there exists an irreducible curve $C_x$ with $D\cdot C_x =0$, then the same conclusion holds for any point of $X$.
In this case, $D$ is semiample by~\cite{lop}*{Theorem~2.7(i)}.
\\
Therefore, we can assume that through the general point of $X$ there is no curve $C$ such that $D\cdot C=0$.
By the structure of the Hilbert scheme, there are at most countably many schemes $T_i$, $i \in I$ whose general points parametrize irreducible curves $C$ with $D\cdot C=0$.
If each $T_i$ is $0$-dimensional then property 2) in the statement of the lemma holds.
Hence, we are left to consider the case that at least one of the $T_i$, say $T_1$, contains a 1-parameter family of curves and the union of all curves parametrized by $T_1$ does not cover $X$. 
In particular, $T_1$ contains a curve $L$ such that the curves parametrized by $L$ span a locus on $X$ of dimension $>1$ so that there is an irreducible surface $S$ covered by curves parametrized by $L$.
Then~\cite{lop}*{Theorem~2.7(i)} implies that $D \cdot S =0$.
On the other hand, we have the following claim.
\begin{claim}
The surface $S$ is an ample divisor in $X$.
\end{claim}
\begin{proof}[Proof of the claim.]
As $X$ does not contain rational curves, $\nef(X) = \Eff(X)$: in fact,
running the $B$-MMP for a big but not nef divisor $B$ would imply the existence of rational curves on $X$, as implied by the Cone Theorem~\cite{kollar-mori}*{Theorem~3.7}.
Therefore, $S$ is an effective nef divisor on $X$ and $S$ is semiample by log abundance for threefolds.
If $S$ is not ample, then $X$ must contain rational curves, cf.~Remark~\ref{rem.rc.fib}, a contradiction to the assumptions of the lemma.
\end{proof}
\noindent
Lemma~\ref{positive.lemma} and the claim imply that $D^2\cdot S>0$, which provides the sought contradiction.
\end{proof}

\section{Proof of Theorem~\ref{main.thm}}\label{sec3}
First let us recall Fujita's vanishing theorem 
from~\cite{fujita-1}*{Theorem~(1)}
and~\cite{fujita-2}*{(5.1) Theorem}, 
see also~\cite{laz}*{Theorem~1.4.35}.

\begin{thm}[Fujita's vanishing]\label{fujita-van}
 Let $X$ be a projective scheme. Let $\mathcal{F}$ be a coherent sheaf and 
 $H$ an ample Cartier divisor on $X$. Then there exists $t_0= t_0(\mathcal{F}, H)$ such that
 \[
 H^i(X, \mathcal{F} \otimes \strutt_X(tH+D))=0
 \]
for all $i>0$, $t \geq t_0$ and any nef Cartier divisor $D$ on $X$.
\end{thm}

We use Fujita's vanishing theorem to derive the following vanishing result for a general coherent sheaf, when we are in the presence of a nef divisor restricting to an ample divisor along an hyperplane section.

\begin{lem}
\label{vanish.lemma}
Let $X$ be a projective scheme. Let $\mathcal{F}$ be a torsion-free  coherent sheaf, $H$ a very ample Cartier divisor and $D$ a nef Cartier divisor on $X$. 
Assume that there exists a very general member $G \in |H|$ of the linear system such that $D|_G$ is ample.
Then there exists a positive integer $r_0$ such that 
\[
H^i(X, \mathcal{F} \otimes \strutt_X(rD)) =0
\]
for all $i\geq 2$ and $r \geq r_0$.
\end{lem}

\begin{proof}
By Fujita's vanishing theorem, we know that there exists $t_0= t_0(\mathcal{F}, H)$ such that
 \begin{align}
 \label{lemma.eq.1}
 H^i(X, \mathcal{F} \otimes \strutt_X(tH+mD))=0
 \end{align}
for all $i>0$, $m\geq 0$ and $t \geq t_0$.
Moreover, we can assume that $D|_{G'}$ is ample for a general
$G' \in |t_0H|$. 
In fact, since $D|_G$ is ample, $D$ is also ample when restricted to the non reduced scheme associated to $t_0G \in |t_0H|$ by~\cite{laz}*{Proposition 1.2.16}. 
Since being ample is an open property along the fibers of a proper morphism, the claim follows.
\\
Let us consider the following short exact sequence 
\[
0 \to \fol \otimes \strutt_X(mD) \to \fol\otimes\strutt_X(mD+G') \to \mathcal{G} \otimes \strutt_{G'}((mD+G')\vert_{G'}) \to 0,
\]
where $\mathcal{G} := \fol / \fol \otimes \strutt_X(-G')$. 
The exactness on the left follows from the torsion-freeness assumption for $\fol$.
The vanishing in~\eqref{lemma.eq.1} implies that 
\[
H^i(G', \mathcal{G} \otimes \strutt_{G'}((mD+G')\vert_{G'})) \simeq H^{i+1}(X, \fol\otimes\strutt_X(mD))
\]
for all $i\geq 1$ and  $m \geq 0$.
Since $D|_{G'}$ is ample, by Serre's vanishing there exists a positive integer $m_0$ such that
\[
H^i(G', \mathcal{G} \otimes \strutt_{G'}((mD+G')\vert_{G'})) =0
\]
for all $i\geq 1$ and $m\geq m_0$, thus concluding the proof of the lemma.
\end{proof}

\begin{rem}
When $X$ is a Calabi--Yau manifold of dimension $n$ 
and $\mathcal{F}$ is its cotangent bundle, 
we can give a quick alternative proof of the vanishing result in Lemma~\ref{vanish.lemma} without referring to Fujita's vanishing theorem.
\\
In fact, let $\pi\colon Y:=\mathbb P(\Omega^1_{X}) \to X$ be the projection morphism  such that $\pi_\ast\mathcal O_{Y}(1)=\Omega^1_{X}$. Note that $\mathcal O_{Y}(1)$ is $\pi$-ample.
Therefore, there exists a fixed integer $t_1$ such that $\mathcal O_{Y}(1)\otimes \mathcal O_{Y}(t\pi^\ast H)$ is ample for all $t\geq t_1$. 
Equivalently, for such $t$, $\mathcal G_t:=\Omega^1_{X}\otimes \mathcal O_X(tH)$ is ample. 
Hence,
\begin{equation}
\label{eq-van}
H^i (X, \omega_X\otimes \mathcal G_t \otimes \det \mathcal G_t \otimes \mathcal O_X(D'))=0
\end{equation}
for all $i\geq 1$ and any nef divisor $D'$ by~\cite{lop}*{Lemma 2.11}.
By assumption, $\det\Omega_X^1\simeq \omega_X\simeq \mathcal O_X$.
Therefore, $ \det \mathcal G_t\simeq \det \Omega^1_{X} \otimes \mathcal O_X(ntH)\simeq \mathcal O_X(ntH)$.
Taking $\bar{t}=(n+1)t_1$, it follows from~\eqref{eq-van} that
\begin{equation*}\label{eq-van-1}
H^i (X, \Omega^1_{X}\otimes\mathcal O_X(\bar{t}H+D'))=
H^i (X, \omega_X\otimes \mathcal G_{t_1} \otimes \det \mathcal G_{t_1} \otimes \mathcal O_X(D'))=0
\end{equation*}
for all $i\geq 1$ and any nef divisor $D'$. 
Taking $D'=D+(t-\bar{t})H$ then completes the argument.
\end{rem}

\begin{proof}[Proof of Theorem~\ref{main.thm}.]
Let us assume by contradiction that $D$ is not semiample.
Then, by~\cite{lop}*{Proposition~2.2(iv)} and Serre duality, 
\begin{equation}
\label{RR.eq1}
\chi(X, \Omega_X^q \otimes \strutt_X(sD)) = (-1)^{q}\frac{c_3(X)}{2}, \quad \text{for} \ q=1, 2, \ \forall s \in \mathbb{Z}.
\end{equation}
The assumption in the statement of the theorem implies that there exists a positive integer $s_0$ such that for any $s \geq s_0$
\begin{align}
\label{eqn.van.proof.1.6}
H^i(X, \Omega_X^q \otimes \strutt_X(sD))=0, \ \text{for} \ i\geq 2, \ q=1, 2,
\end{align}
by Lemma~\ref{vanish.lemma}.
Hence, Proposition~\ref{LP.van.prop} and~\eqref{eqn.van.proof.1.6} imply that for any $s \geq s_0$, we can rewrite~\eqref{RR.eq1} in the following ways
 \begin{align*}
 -\frac{c_3(X)}{2}&=\chi(X, \Omega_X^1 \otimes \strutt_X(sD)) = -h^1(X, \Omega_X^1 \otimes \strutt_X(sD))\leq 0,\\
  \frac{c_3(X)}{2}&=\chi(X, \Omega_X^2 \otimes \strutt_X(sD)) = -h^1(X, \Omega_X^2 \otimes \strutt_X(sD))\leq 0.
 \end{align*}
 These two inequalities together imply that $c_3(X) = 0$, which gives the desired contradiction.
\end{proof}

\begin{rem}
Let $X$ be a Calabi--Yau threefold and let $\mathcal{F}$ be the cotangent bundle of $X$. 
The result of Theorem~\ref{main.thm} was implicitly claimed by Wilson in the proof of~\cite{wilson}*{Theorem~2.3}. 
As the proof of Wilson's claim is sketched loosely, Lazi\'c--Oguiso--Peternell tried to reprove this result in~\cite{lop}, but ended up needing to add the extra assumption (i) in~\cite{lop}*{Proposition 5.2} to obtain the desired statement. 
Our result shows that this assumption is indeed not necessary. 
\end{rem}

\section{Proof of Theorem~\ref{lop.thm}}
\label{sect.proof.thm.lop}

In this section, we aim to prove Theorem~\ref{lop.thm}.
The variety $X$ will be a Calabi--Yau threefold and $D$ will be a nef divisor of numerical dimension $2$ on $X$.
Let us recall that a prime divisor $S \subset X$ is orthogonal to $D$ if $S \cdot D=0$, or, equivalently, $D \vert_S \equiv 0$.

A birational morphism $\phi\colon X \to Z$ is said to be {\it a good Calabi--Yau model} for a nef divisor $D$ if $Z$ is a $\mathbb Q$-factorial $K$-trivial threefold with canonical singularities and the exceptional locus of $\phi$ contains all and only the divisors orthogonal to $D$, cf.~\cite{lop}*{Definition~8.2}. 

In~\cite{lop}*{Theorem~8.3}, it is shown that there exists a smooth Calabi--Yau threefold $X'$ isomorphic in codimension one to $X$ on which the strict transform $D'$ of $D$ is nef, and a good Calabi--Yau model $\phi' \colon X' \to Z$ for $D'$.
In particular, it is immediate to see that $D$ and $D'$ are crepant birational to each other, and $D$ is semiample if and only if $D'$ is;
furthermore, a prime divisor $S \subset X$ is orthogonal to $D$ if and only if its strict transform $S' \subset X'$ is orthogonal to $D'$.

As $X$ and $X'$ are isomorphic in codimension one, their Hodge numbers coincide, cf.~\cite{bat2}, so that $c_3(X)=c_3(X')$.
Moreover, any prime divisor $T \subset X$ is birational to its strict transform $T' \subset X'$.
Letting $S_1, \dots, S_r$ be the prime divisors orthogonal to $D$ and $S'_1, \dots, S'_r$ be the strict transforms of the $S_i$ on $X'$, then the irregularity $g_j$ of a resolution of $S_i$ coincide with that of a resolution of $S'_i$.

Hence, in view of these observations, in order to prove Theorem~\ref{lop.thm}, we can substitute $X$ (resp. $D$) with $X'$ (resp. $D'$) and assume that $X$ is endowed with a good Calabi--Yau model for $D$.
Thus, it suffices to prove the following theorem.

\begin{thm}
\label{lop.thm1}
Let $X$ be a Calabi--Yau threefold and $D$ a nef divisor on $X$ with $\nu(D)=2$, such that there exists a birational morphism $\phi\colon X \to Z$ which is a good Calabi--Yau model for $D$. 
Let $S_1, \ldots, S_r$ be all the prime divisors on $X$ orthogonal to $D$. 
Let $g_j$ be the irregularity of a resolution of $S_j$. If
\[
\frac{c_3(X)}{2} \neq r-\sum_{j=1}^r g_j,
\]
then $D$ is semiample.
\end{thm}

In order to prove Theorem~\ref{lop.thm1}, we will need the following generalization of~\cite{lop}*{Proposition~8.13}.

\begin{lem}
\label{lem.various.vans}
Notation as in Theorem~\ref{lop.thm1}. Fix a resolution $\tau\colon Y \to X $ such that the exceptional set of the morphism $\pi\colon Y \to Z$ where $\pi:=\phi \circ \tau$, is a simple normal crossings divisor $E = \sum_{j=1}^s  E_j$ on $Y$. Suppose that $D$ is not semiample.
Then,
\[
H^q\big(Y, \Omega^1_Y(\log E) \otimes \pi^*\mathcal{O}_Z(mD_Z)\big) = 0, \ \forall q  \ \text{and} \ \forall m \gg 0.
\]
\end{lem}

\begin{proof}
As showed in~\cite{lop}*{Proposition 8.13}, it suffices to prove the following vanishings:
\begin{enumerate}
\item[(i)] $H^1\big(Z, \Omega^{[1]}_Z \otimes \mathcal{O}_Z(mD_Z)\big) = H^3\big(Z, \Omega^{[1]}_Z \otimes \mathcal{O}_Z(mD_Z)\big)=0$ for any $m\ll0$,
\item[(ii)] $\chi\big(Z,\Omega^{[1]}_Z \otimes \mathcal{O}_Z(mD_Z)\big) =  0 $ for all integers $m$,
\item[(iii)] $H^q\big(Z,\Omega^{[1]}_Z \otimes \mathcal{O}_Z(mD_Z)\big) = 0$ for all $q$ and all integers $m$ such that $| m | \gg 0$, and
\item[(iv)] $ R^i\pi_\ast\Omega^1_Y(\log E) = 0$ for $i=1,2$.
\end{enumerate}
We only sketch the proof of (i) and (ii), since the proof of the other two is exactly the same as that in~\cite{lop}*{Proposition 8.13} and does not use the extra hypotheses in the statement of the proposition. 
Let us recall that 
\begin{align}
\label{eqn.gkkp}
\pi_\ast\Omega^p_Y(\log E) \simeq 
\Omega^{[p]}_Z\quad\text{and}\quad \phi_\ast\Omega^p_X \simeq 
\Omega^{[p]}_Z\quad \forall \ 0\leq p\leq 3,
\end{align} 
as in~\cite{lop}*{Equation (37)}.

\begin{enumerate}
    
\item[(i)]
As $K_Z \sim 0$, by Serre duality, it follows that 
\[
H^1\big(Z, \Omega^{[1]}_Z \otimes \mathcal{O}_Z(mD_Z)\big) \simeq
H^2\big(Z, \Omega^{[2]}_Z \otimes \mathcal{O}_Z(-mD_Z)\big),
\quad \forall m \in \mathbb{Z}.
\]
Proposition~\ref{prop.ample.surf} and Lemma~\ref{vanish.lemma} imply that
\begin{align}
\label{eqn.van.h2}
H^2(\Omega^{[2]}_Z \otimes\mathcal{O}_Z(m D_Z)) =0, \quad 
\forall m \gg 0.
\end{align}
Hence, the first part of (i) follows by combining these two observations.
\\
The second vanishing claimed in (i) follows instead from the following isomorphisms
\begin{align*}
&H^3\big(Z, \Omega^{[1]}_Z \otimes \mathcal{O}_Z(mD_Z)\big)\\
\simeq & {\rm Hom}\big(\Omega^{[1]}_Z \otimes
\mathcal{O}_Z(mD_Z), \mathcal{O}_Z\big) & \hfill[\text{by Serre duality}]\\
\simeq & H^0\big(Z,\Omega^{[2]}_Z \otimes \mathcal{O}_Z(-mD_Z)\big)
&  \hfill[\text{since} \ K_Z \sim 0]\\
\simeq & H^0\big(X,\Omega^2_X \otimes \mathcal{O}_X(-mD)\big)=0 & \hfill [\text{by} \ \eqref{eqn.gkkp}]\\
=& 0. & \hfill [\text{by Proposition} \ \ref{LP.van.prop}]
\end{align*}
    
    \item[(ii)]
By the same proof of~\cite{lop}*{Proposition 8.13, (iii)}, we obtain that
\[
\chi\big(Z, \Omega^{[1]}_Z \otimes\mathcal{O}_Z(m D_Z)\big) \ \text{is independent of } m.
\]
By Proposition~\ref{LP.van.prop} and~\eqref{eqn.gkkp}, $H^0(Z, \Omega^{[1]}_Z \otimes\mathcal{O}_Z(m D_Z))= 0$, for all $\vert m \vert \gg 0$.
Moreover, by Proposition~\ref{prop.ample.surf} and Lemma~\ref{vanish.lemma},
\begin{align}
\label{eqn.van.fuj.appl}
H^i(\Omega^{[j]}_Z \otimes\mathcal{O}_Z(m D_Z)) =0, \quad 
\forall m \gg 0, \ i=2,3, \ j=1,2.
\end{align}
Hence,~\eqref{eqn.van.fuj.appl} and Serre duality, together with (i), imply that
\[
\chi\big(Z, \Omega^{[1]}_Z \otimes\mathcal{O}_Z(m D_Z)\big) 
\begin{cases}
\leq 0 & \text{for } m \gg 0,\\
\geq 0 & \text{for } m \ll 0, 
\end{cases}
\]
which concludes the proof.
\end{enumerate}
\end{proof}

\begin{proof}[Proof of Theorem~\ref{lop.thm}]
With the notation introduced above, the proof of~\cite{lop}*{Theorem~8.5} applies almost verbatim after replacing~\cite{lop}*{Theorem~8.13} with our Lemma~\ref{lem.various.vans}. 
\end{proof}

In the proof of Lemma~\ref{lem.various.vans}, we used the following result.

\begin{prop}
\label{prop.ample.surf}
Let $X$ be a Calabi--Yau threefold.
Let $D$ be a nef divisor of numerical dimension $2$ on $X$ which is not semiample.
Assume that there exists a good Calabi--Yau model $\phi \colon X \to Z$ for $D$.
Let $D_Z$ be a $\mathbb Q$-Cartier divisor on $Z$ such that $D \sim_\mathbb{Q} \phi^\ast D_Z$.
Let $H_Z$ be a very ample divisor on $Z$. Then
for a very general $G \in \vert H_Z \vert$, $D_Z\vert_G$ is ample.
\end{prop}

\begin{proof}
We give a proof by contradiction.\\
Assume that for all $G \in \vert H_Z \vert$, $D_Z \vert_G$ is not ample.
By Lemma~\ref{positive.lemma}, $D_Z^2\cdot G>0$.
Hence, by the Nakai--Moishezon criterion, for any $G \in \vert H_Z \vert$ there exists a curve $C_G \subset G$ such that $D_Z \cdot C_G=0$.
In particular, as the Hilbert scheme of $Z$ has countably many components, there exists a family $\mathcal C$ of curves that have intersection $0$ with $D_Z$ on $Z$, such that the dimension of $\mathcal C$ is $\geq 2$.
Let
\[
S_\mathcal{C}= \bigcup_{C \in \mathcal C} C \subset Z.
\]
Then $\dim S_\mathcal{C} \geq 2$.
\\
If $S_\mathcal{C} = Z$, we show that $D$ is semiample, which prompts the desired contradiction.
Indeed, by taking the strict transform $C_X$ of the general curve $C \in \mathcal C$ on $X$, it follows that $X$ itself is covered by curves having intersection $0$ with $D$ and we can apply then~\cite{lop}*{Theorem~2.7(ii)} to prove that $D$ is semiample.
The strict transform $C_X$ of the general element $C \in \mathcal C$ exists, since $\mathcal C$ covers $Z$ and so $C_X$ is not contained in the exceptional locus of $\phi$.
\\
If $S_\mathcal{C} \subsetneq Z$, let $S'$ be an irreducible component of $S_\mathcal{C}$.
Then, it follows from~\cite{lop}*{Theorem~2.7(i)} that \begin{align}
\label{eqn.surf.inters}
S'_X \cdot D = 0, 
\quad \text{where} \ S'_X:=\phi^{-1}_\ast S'.
\end{align}
This prompts the desired contradiction, since~\eqref{eqn.surf.inters} implies that $S'_X$ would have to be contracted by $\phi$, by definition of a good Calabi--Yau model.
\end{proof}

\section{Strictly nef divisors}
\label{sec5}

Let us recall that a strictly nef divisor $D$ on a normal variety $X$ is a nef Cartier divisor such that $D\cdot C>0$ for any curve $C$ on $X$.
In general, a strictly nef divisor $D$ is not necessarily ample, as shown
by many classical example, e.g., the examples of Mumford and Ramanujam, see~\cite{hart-1}*{appendix to Chapter I}.
On the other hand, given a strictly nef divisor $D$ on a projective variety $X$, the ampleness of $D$ is implied by $D$ being semiample, as then strict nefness of $D$ forces the Iitaka fibration $\pi\colon X\to Z$ of $D$ to be finite.

While, as we have just discussed, we cannot expect strictly nef divisors to always be ample, the situation is expected to improve if we allow twisting with the canonical bundle, as it is often the case when working with positivity properties in algebraic geometry.
In this context, Serrano~\cite{serrano} proposed the following conjecture that aims to bound the distance between ample and strictly nef divisors once we allow such twisting.
This conjecture can be also thought as a weak analogue of Fujita's Conjectures, cf.~\cite{laz}*{\S~10.4.A}, for strictly nef divisors.

\begin{conj}\label{conj}
Let $X$ be a projective manifold of dimension $n$ and $D$ a strictly nef divisor on $X$. Then $K_X+tD$ is ample for any real number $t>n+1$.
\end{conj}

In~\cite{serrano}, Serrano showed that Conjecture~\ref{conj} holds in dimension 2, and that in dimension 3 the only unknown cases are the following two:
\begin{itemize}
\item[(i)] $X$ is a Calabi--Yau threefold and $c_2(X)\cdot D=0$, and
\item[(ii)] $\kappa(K_X)=-\infty$, and either the irregularity $q(X)\leq 1$ or else $q(X)=2$ and $\chi(\strutt_X)=0$.
\end{itemize} 
\noindent
A decade later, Campana--Chen--Peternell~\cite{ccp} ruled out the latter case and proved some partial results for projective manifolds of dimension higher than $3$. 
When $X$ is a Calabi--Yau threefold and $D$ is a strictly nef divisor,
as in~\S~\ref{cy.sect}, we can assume that $\kappa(D)=-\infty$
and $c_2(X) \cdot D=0$, which in turn implies that $D^3=0$ and $1\leq \nu(D)\leq 2$.

The results proven in the previous sections allow us to prove that Serrano's conjecture holds on a Calabi--Yau threefold if we assume that the numerical dimension $\nu(D)$ of the strictly nef divisor $D$ is $\neq 1$, when $c_3(X) \neq 0$.

\begin{thm}
\label{sn-thm}
Let $X$ be a Calabi--Yau threefold with $c_3(X)\neq 0$. 
Let $D$ be a strictly nef divisor on $X$ with $\nu(D)\neq 1$.
Then $D$ is ample.
\end{thm}

\begin{proof}
The result follows at once from Corollary~\ref{exist-cor} and Theorem~\ref{main.thm}.
\end{proof}

In the case $\nu(D)=1$, instead, it is possible to prove the following result, valid for Calabi--Yau varieties of any dimension $\geq 3$, which is a special case of~\cite{lp}*{Theorem~6.5}.

\begin{thm}
Let $X$ be a Calabi--Yau manifold of dimension $n\geq 3$
and $D$ a nef divisor on $X$ with $\nu(D)=1$.
Assume that there is a singular metric $h$ on $\strutt_X(D)$
with semipositive curvature current such that the multiplier ideal sheaf $\mathcal I(h) \subsetneq \strutt_X$. 
Then $\kappa(D)\geq 0$.
\end{thm}

Note that for a pseudoeffective Cartier divisor $D$, thus including the case where $D$ is nef, we can always find a singular metric $h$ on $\strutt_X(D)$ with semipositive curvature current. 
In general, it is rather hard to determine whether the multiplier ideal sheaf $\mathcal I(h)$ of the metric $h$ is contained in $\strutt_X$ strictly or not.
\begin{cor}
Let $X$ be a Calabi--Yau threefold
and $D$ a strictly nef divisor on $X$.
If there is a singular metric $h$ on $\strutt_X(D)$
 with semipositive curvature current such that 
$\mathcal I(h) \subsetneq \strutt_X$, then $\nu(D)\neq 1$.
\end{cor}

We conclude this section with a final remark on the case  $\nu(D)=1$.

\begin{rem}
Let $X$ be a Calabi--Yau threefold
and $D$ a strictly nef divisor on $X$ with $\nu(D)=1$.
Then $X$ cannot contain any del Pezzo, $K3$ or abelian surfaces. 
\\
Indeed, if $S \subset X$ is a del Pezzo surface, then $D|_S$ is semiample by the basepoint-free theorem.
As $D|_S$ must also be strictly nef, $D|_S$ is ample which contradicts the assumption that $\nu(D)=1$.
It is well-known that, on abelian and $K3$ surfaces strictly nef divisors are semiample, hence ample;
thus, the same argument as in the del Pezzo case shows that such types of surfaces cannot be contained in $X$.

As a consequence of this simple observation, the existence of a strictly nef divisor $D$ with $\nu(D)=1$ implies that there cannot possibly exist morphisms $\phi \colon X \to Y$ such that:
\begin{itemize}
    \item 
    $\phi$ is a birational morphism contracting a del Pezzo surface to a point; or
    
    \item 
    $\phi$ is a morphism whose generic fibre is an abelian or $K3$ surface.    
\end{itemize}

In particular, there is no nef effective divisor $E$ on $X$ such that $\nu(E)=1$: 
indeed, $E$ would then be semiample and its Iitaka fibration would give a fibration in abelian/K3 surfaces, see~\S~\ref{cy.sect}.
\\
Conversely, if there is a contraction $\phi$ of $X$ 
as above, then there is no strictly nef divisor of numerical dimension 1 on $X$.
This type of observation will also be used in the analysis of the existence of rational curves in the next section.
\end{rem}

\section{Existence of rational curves}\label{sec6}

Another classical problem in the study of K-trivial varieties is to determine whether or not they contain rational curves.
Existence of rational curves is fully established only in 
dimension $2$ thanks to Bogomolov and Mumford~\cite{mm}, while only partial results are known in dimension $3$ and higher, as we discuss below.
Even though it is widely believed that every Calabi--Yau manifold should contain rational curves, already in the threefold case this appears to 
be a very difficult problem.

In his study of the birational structure of Calabi--Yau threefolds, 
Oguiso proposed Conjecture~\ref{og.rat.curv.conj.intro} which is a weakened version of the classical question on the existence of rational curves.
We recall it here for the reader's convenience.

\begin{conj}
\label{og.rat.curv.conj}
	Let $X$ be a Calabi--Yau threefold. 
	Assume that there exists a non-trivial Cartier divisor 
	$D$ contained in the boundary of $\nef(X)$.
	Then $X$ contains a rational curve.
\end{conj}

The existence of a divisor $D$ satisfying the hypotheses of Conjecture~\ref{og.rat.curv.conj} is a rather strong assumption. 
First of all, it requires that $h^2(X) >1$. 
It also implies that the boundary of the nef cone contains rational points other than the origin. 
The lack of rational points in the boundary of the nef cone is a rather typical situation if we do not assume that $\pi_1(X)=\{1\}$ or  $h^i(X, \mathcal{O}_X)=0$ for $0< i<\dim X$.
For example, the nef cone of an abelian variety contains non-zero rational points if and only if it admits a morphism to a lower dimensional abelian variety; 
in~\cite{oguiso-1}, Oguiso has constructed examples of hyperk\"ahler manifolds and of Calabi--Yau threefolds of Picard number 2 whose nef cones do not contain rational points other than zero in the boundary. 
Nonetheless, abelian varieties do not contain rational curves, while Verbitsky proved in~\cite{ver2} that no hyperk\"ahler manifold is Kobayashi hyperbolic.
Moreover, the analogous of Conjecture~\ref{og.rat.curv.conj} is expected to hold also for hyperk\"ahler manifolds.

In a fundamental series of works~\cites{wilson1,wilson2, wilson2.err}, Wilson studied the structure of the nef cone of a Calabi--Yau threefold. 
In particular, he showed that if $\rho(X)$ is sufficiently big ($\rho(X) \geq 13$), then there always exists a birational morphism $f \colon X \to Y$. 
In particular, this forces the existence of a non-zero rational point in the boundary of the nef cone of $X$.
Moreover, as already mentioned in Remark~\ref{rem.rc.fib}, the exceptional locus of $f$ is uniruled, so that rational curves exist in $X$. 
Thus, to solve Conjecture~\ref{og.rat.curv.conj}, we just need to focus on Calabi--Yau varieties of low Picard rank, that is, of small second Betti number.
Building on the seminal works of Wilson, Diverio--Ferretti proved the following result in~\cite{df}*{Theorem~1.2} which partially solves Conjecture~\ref{og.rat.curv.conj}.

\begin{thm}\label{div.ferr.thm}
	Let $X$ be a Calabi--Yau threefold. 
	Assume that there exists a non-trivial Cartier divisor $D$ 
	contained in the boundary of  $\nef(X)$ and $h^2(X) \geq 5$.
	Then $X$ contains a rational curve.
\end{thm}

We briefly explain some of the arguments contained in the proof of the above theorem.

By contradiction, if there is no rational curve on $X$, we can assume that $\nef(X)=\Eff(X)$ and $D^3 =0=c_2(X) \cdot D$, cf. the proof of  Lemma~\ref{h.sch.lemma} and \S~\ref{cy.sect};
thus, we can assume that the cohomology class of $D$ lies on the intersection between the cubic hypersurface
\begin{align}
\label{def.W}
W := \{T \in H^2(X, \mathbb{R}) \; | \; T^3=0\} \subset 
	H^2(X, \mathbb{R})
\end{align}
and the hyperplane
\begin{align}
    \label{eqn.def.perp}
c_2(X)^\perp :=
\{  H \in H^2(X, \mathbb{R}) \; | \; c_2(X) \cdot H = 0\}.
\end{align}
In this case, the Hirzebruch--Riemann--Roch formula in~\eqref{RR.cy.eq.div} implies that 
\[
\chi(X, \mathcal{O}_X(mD))=0, \quad \forall m \in \mathbb{Z}.
\] 
It is then not clear in general how to produce sections in any of the linear systems $|mD|$. 
Let $\widetilde{W}\subset \mathbb{P}(H^2(X, \mathbb{R}))$ be the projectivization of $W$. If $D^2=0$, then the point $p_D$ associated to $D$ in $\mathbb{P}(H^2(X, \mathbb{R}))$ is a singular point of $\widetilde{W}$ and one can show that $\widetilde{W}$ contains a dense set of rational points. 
In particular, as $\nef(X) \subset W$, this implies that $X$ carries an elliptic fibration, and, thus, by Remark~\ref{rem.rc.fib} it must contain a rational curve.

When $\nu(D)=2$, in~\cite{df}, the authors carefully analyze the configuration of $\widetilde{W}$ and $c_2(X)^\perp$, together with results of Wilson, to conclude the proof of Theorem~\ref{div.ferr.thm}.

We show that when $c_3(X) \neq 0$ it is actually possible to fully prove Oguiso's conjecture.
We use the results of~\S~\ref{sec3} to show that under these assumptions the divisor $D$ contained in the boundary of the nef cone will be semiample.
Hence, by Remark~\ref{rem.rc.fib}, $X$ will then contain a rational curve.

\begin{thm}\label{rc-thm}
Let $X$ be a Calabi--Yau threefold. 
Assume that there exists a non-trivial Cartier divisor $D$ contained in the boundary of  $\nef(X)$.
Then, $X$ contains a rational curve provided that either $c_3(X) \neq 0$ or the second Betti number $h^2(X) \neq 2,3,4$.
\end{thm}

\begin{proof}
{\bf Case $c_3(X) \neq 0$}.
By the above discussion, we can assume that $c_2(X)\cdot D=D^3 = 0 \neq D^2$.
If $X$ does not contain any rational curve, then by Lemma~\ref{h.sch.lemma} the set $\mathcal{C}_D$ is at most countable. 
Hence, Lemma~\ref{curve.lemma} and Theorem~\ref{main.thm} imply that $D$ is semiample, 
which prompts a contradiction by Remark~\ref{rem.rc.fib}.
\\
{\bf Case $c_3(X) = 0$}. 
The result follows from Theorem~\ref{div.ferr.thm}.
\end{proof}

The material in the next section contains some considerations about possible strategies to prove the missing cases.

\section{Kawamata--Morrison Conjecture}\label{sec7}

In the final section we explain how Conjecture~\ref{og.rat.curv.conj} is closely intertwined with another important conjecture regarding Calabi--Yau manifolds, the so-called Kawamata--Morrison Cone Conjecture,~\cites{mor.km,kaw.km}.
\begin{conj}[Kawamata--Morrison Cone Conjecture]\label{kaw.morr.conj}
Let $X$ be a Calabi--Yau manifold and $\Aut(X)$ the group of automorphisms of $X$. There exists a rational polyhedral cone $\Pi$ which is a fundamental domain for the action of ${\rm Aut}(X)$ on ${\rm Nef}^e(X)$, in the sense that
    \begin{itemize}
        \item [a)] ${\rm Nef}^e(X)=\bigcup_{g \in {\rm Aut}(X)} g^\ast  \Pi$; and
        \item[b)] ${\rm Int}\Pi \cap  {\rm Int}g^\ast\Pi = \emptyset$ unless $g^\ast ={\rm Id}_{{\rm N}^1(X)}$.
    \end{itemize}
Moreover, the number of ${\rm Aut}(X)$-equivalence classes of faces of the cone ${\rm Nef}^e(X/Y)$ corresponding to birational contractions or fiber space structures is finite.    
\end{conj}

There is an analogous version of the conjecture involving the birational automorphism groups and its action on the cone of movable effective divisor of $X$.
There is now a fairly extensive literature on the Cone conjecture;
for the case of Calabi--Yau threefolds that we treat here, the most relevant papers are~\cites{lp13, lop1}.

The above conjecture predicts that if $\Aut(X)$ is finite then the cone on which that group acts should be rational polyhedral. 
It follows that, if $\nef(X)$ is polyhedral, then those divisors contained in the facets of $\nef(X)$ not contained in the hyperplane $c_2(X)^\perp$ are semiample. 
In that case, $X$ will contain rational curves as in \S~\ref{cy.sect}.
Moreover, assuming Conjecture~\ref{kaw.morr.conj} holds, Conjecture~\ref{og.rat.curv.conj} holds in some of the 
cases that are not covered by Theorem~\ref{rc-thm}, as we proceed to explain now.

We follow the notation introduced in the previous section. 

\begin{prop}
	Let $X$ be a Calabi--Yau threefold.
	Assume that Conjecture~\ref{kaw.morr.conj} holds true.
	Then Conjecture~\ref{og.rat.curv.conj} holds true unless $\rho(X) = 3,4$ and $\widetilde{W}$ is the union of a hyperplane $P$ and a quadric $Q$ of rank $\rho(X)$ such that $P=c_2(X)^\perp$, and $Q$ does not contain any rational points.
\end{prop}

\begin{proof}
We may assume, by contradiction, that $X$ does not contain any rational curves and proceed to a case by case analysis.
Thus, $\nef(X)= \Eff(X)$ and the boundary of the nef cone is contained in $W$ and contains no semiample divisor.
Let $D$ be a non-zero Cartier divisor contained in the boundary of $\nef(X)$; thus, $D^3=0=c_2(X)\cdot D$.
By abusing notation, we will use $D$ to denote also the rational point in $\mathbb{P}(H^2(X, \mathbb{R}))$ associated to it, so that $D \in \widetilde W$.
\\
{\bf Case $h^2(X)=\rho(X)=2$}. Then the hypersurface $W \subset H^2(X, \mathbb R)$ is a union of three lines, $l_1, l_2, l_3$;
one of the lines, say $l_1$ is generated by the class of $D$, as $D^3=0$; hence $c_2(X)^\perp=l_1$, while the other two lines, $l_2, l_3$, are defined over a quadratic extension of $\mathbb{Q}$ and are conjugated under the $\mathbb{Z}_2$-action of Galois group. 
As the matrix $M_g$ corresponding to the action of an element of $g \in \Aut(X)$ on ${\rm N}^1(X)_\mathbb{R}$ is defined over the underlying integral structure, and $l_1$ is an eigenspace of eigenvalue $1$, then $M_g$ must be the identity as it has to preserve $\nef(X)$.
Hence, Conjecture~\ref{kaw.morr.conj} implies that $\nef(X)$ is polyhedral, which then implies $\partial \nef(X) \not \subset W$.
\\
{\bf Case $h^2(X)=\rho(X)=3$}.
Then the possible configurations for the pair $(\widetilde W, c_2(X)^\perp)$, excluding the one in the statement of the proposition, are the following:
\begin{enumerate}
    \item 
    $\widetilde W$ is the union of three distinct lines $H_1, H_2, H_3$, while $c_2(X)^\perp$ is a line;
    
    \item 
    $\widetilde W$ is the union of a line $P$ and a conic $C$, and $c_2(X)^\perp$ is a line. If the rank of $C$ is maximal, then $C$ contains rational points;
    
    \item 
    $\widetilde W$ is geometrically irreducible and $c_2(X)^\perp$ is a line $L$.
\end{enumerate}
We will show that in all these cases, $\nef(X)$ must contain a semiample divisor.
In case (1) either the three
hyperplanes $H_1, H_2, H_3$, are all defined over $\mathbb{Q}$ and so the boundary of $\nef(X)$ contains a nef divisor $E$ with $c_2(X)\cdot E>0$; 
or, up to relabeling the $H_i$, $H_1$ is defined over $\mathbb{Q}$ and $H_2, H_3$ are defined over a real quadratic extension of $\mathbb{Q}$ and are conjugated under the $\mathbb{Z}_2$-action of the Galois group. 
In this case, either $H_1 \neq c_2(X)^\perp$, so that, again, $H_1$ will contain a nef and semiample divisor, or the point $p_F=H_2 \cap H_3$ is defined over $\mathbb{Q}$ and it corresponds to  projectivization of a nef divisor $F$ with $c_2(X) \cdot F >0$, and again we can conclude as above.
\\
In case (2), the boundary of the nef cone is contained in $\widetilde W$, hence either $P$ is tangent to $C$ and the nef cone is contained in the oval given by the interior of $C$, or the nef cone is contained in one of the subsets in which the line $P$ partitions the oval given by the interior of $C$.
In the former case, then $P \cap C$ is a rational point, hence $C$ is rational over $\mathbb Q$ and we are done.
In the latter case, if $c_2(X)^\perp \neq P$ then $P$ will contain rational points corresponding to semiample divisors.
Hence, we can assume that $c_2(X)^\perp = P$.
Moreover, if $C$ has maximal rank and contains a rational point, then we can conclude as above.
While, if rank $C=2$ then $C$ is a cone and its vertex corresponds to a nef divisor $D'$ with $D^2=0\neq c_2(X)\cdot D$. 
Hence $D'$ is semimaple.
\\
In case (3), as $\widetilde{W}$ is geometrically irreducible, 
it follows immediately that the action of $\Aut(X)$ on the second cohomology embeds in the automorphisms
of the pair $(\widetilde{W}, c_2(X)^\perp_{|\widetilde{W}})$, that are finite, as $\widetilde W$ is a cubic hypersurface and we are furthermore fixing the linear section $c_2(X)^\perp$;
hence, the finiteness of the representation of the automorphisms of $X$ on ${\rm N}^1(X)_\mathbb R$ implies that the nef cone of $X$ is rational polyhedral and in particular there are some nef non-ample divisors on $X$ that are semiample.
\\
{\bf Case $h^2(X)=\rho(X)=4$}.
Then,~\cite{wil.aq}*{Proposition~2.1} implies that if we assume that $c_2(X)$ is not $>0$ on $\nef(X) \setminus \{0\}$ and that $\nef(X)$ does not contain any semiample divisors -- which is the case here, in view of \S~\ref{cy.sect} -- then $\widetilde{W}=Q \cup P \subset \mathbb{P}^3_\mathbb{R}$, where $Q$ is a rank 4 quadric with no rational points and $P=c_2(X)^\perp$.
\end{proof}

To conclude, let us explain how Conjecture~\ref{kaw.morr.conj} is also connected to the question of whether or not Calabi--Yau threefolds are Kobayashi hyperbolic.

By~\cite{brody}, a compact complex space is Kobayashi hyperbolic if
and only if it does not contain any non-constant holomorphic curve 
$f\colon \mathbb C\to X$.
A now classical conjecture of Kobayashi's predicts that if a smooth projective variety is Kobayashi hyperbolic then the canonical class $K_X$ is ample on $X$.
In particular, according to this conjecture, no Calabi--Yau manifold is Kobayashi hyperbolic.
Kobayashi's conjecture is known to hold for smooth projective surfaces.
In dimension $3$, the only missing case for the verification of the conjecture is that of Calabi--Yau threefold, see~\cite{div} for more details about the state of the art.

It is not hard to see that if Conjecture~\ref{og.rat.curv.conj} holds then also Kobayashi's conjecture holds in dimension 3.
Even though we are not able to show that if Conjecture~\ref{kaw.morr.conj} holds then also Conjecture~\ref{og.rat.curv.conj} holds, it is still possible to show that $X$ is indeed not Kobayashi hyperbolic at least when $\rho(X)\geq 2$.

\begin{prop}
\cite{div}*{Proposition}
	Assume that Conjecture~\ref{kaw.morr.conj} holds.
	Then the Kobayashi conjecture holds true in dimension 3, 
	except possibly if there exists a Calabi--Yau threefold of 
	Picard number one 
	which is hyperbolic.
\end{prop}


\end{document}